%% file: ms.tex
\documentclass{ifacconf}

\usepackage{graphicx}      
\usepackage{natbib}        
\usepackage{amsmath}	
\usepackage{epstopdf}
\usepackage{fancyhdr}
\input{macros/symbols.tex}
\begin{document}

\lhead{\tiny{NOTICE: this is the author’s version of a work that was accepted for publication in Proceedings of the 15th \textbf{IFAC Symposium on Control in Transportation Systems (CTS 2018)}, 6-8 June, Savona, 
Italy. Changes resulting from the publishing process, such as peer review, editing, corrections, structural formatting, and other quality control mechanisms may not be reflected in this document. Changes may have been made to this work since it was submitted for publication.}}
\chead{}
\rhead{}
\lfoot{}
\cfoot{}
\rfoot{}

\begin{frontmatter}

\title{A novel model-based heuristic for energy-optimal motion planning for automated driving} 


\author[First]{Zlatan Ajanovi{\'c}} 
\author[First]{Michael Stolz} 
\author[Second]{Martin Horn}

\address[First]{Virtual Vehicle Research Center, Department Electrics/Electronics and Software, Inffeldgasse 21/A, 8010 Graz, Austria (e-mail: \{zlatan.ajanovic, michael.stolz\}@v2c2.at).}
\address[Second]{Graz University of Technology, Institute of Automation and Control,
Inffeldgasse 21/B, 8010 Graz, Austria, (e-mail: martin.horn@tugraz.at).}

\begin{abstract}                
Predictive motion planning is the key to achieve energy-efficient driving, which is one of the main benefits of automated driving. Researchers have been studying the planning of velocity trajectories, a simpler form of motion planning, for over a decade now and many different methods are available.  Dynamic programming has shown to be the most common choice due to its numerical background and ability to include nonlinear constraints and models.  Although planning of an optimal trajectory is done in a systematic way, dynamic programming does not use any knowledge about the considered problem to guide the exploration and therefore explores all possible trajectories. 

A* is a search algorithm which enables using knowledge about the problem to guide the exploration to the most promising solutions first. Knowledge has to be represented in a form of a heuristic function, which gives an optimistic estimate of cost for transitioning to the final state, which is not a straightforward task. This paper presents a novel heuristics incorporating air drag and auxiliary power as well as operational costs of the vehicle, besides kinetic and potential energy and rolling resistance known in the literature. Furthermore, optimal cruising velocity, which depends on vehicle aerodynamic properties and auxiliary power, is derived. Results are compared for different variants of heuristic functions and dynamic programming as well.
\end{abstract}


\end{frontmatter}
\thispagestyle{fancy}
\section{Introduction}
Increasing environmental awareness, strict regulations on greenhouse gas emissions and constant desire to increase the range of electric vehicles as well as the big economic benefits motivated increased research interest in the field of energy-efficient driving. So far, many different approaches addressing this topic exist. In \cite{bingham2012impact} authors presented a study which shows that the driving behavior has a rather big influence on energy consumption. It was shown that energy consumption may vary in a range of approx. 30\% depending on driving behavior. Knowledge about the upcoming driving route, the road slope profile and the ability to control the vehicle's propulsion enables the optimization of the velocity trajectory of the vehicle with respect to the energy consumption. Discrete dynamic programming (DP) has been used for over a decade now for finding energy efficient velocity trajectories (\cite{hellstrom2005explicit, hellstrom2010look}). Because of high computational requirements in \cite{ozatay2014cloud} authors used a cloud to compute energy efficient velocity trajectories. By using model predictive control (MPC) to drive vehicles on free roads with up and down slopes, notable fuel savings are shown in  \cite{kamal2011ecological}. MPC was also used to control a hybrid vehicle driving over a hill and performing vehicle following in \cite{vajedi2016ecological}. In \cite{myAMAA2017} authors presented the use of dynamic programming in an MPC-like framework. In  \cite{sciarretta2015optimal} an overview of the existing approaches and the current state-of-the-art can be found. 

One of the first uses of A*-search for autonomous vehicle motion planning was presented in \cite{fraichard1993dynamic}, the result of a Prometheus project. The author used A*-search to find the shortest time motion in the presence of dynamic obstacles and introduced state-time space for dealing with dynamic obstacles. A*-search was also used in DARPA Urban challenge by many teams. Here, especially the application from Stanford has to be mentioned where the authors introduced Hybrid A* search to deal with rounding errors for finding the shortest path (\cite{dolgov2008practical}). Another recent application of A*-search for planning safe trajectories is shown in \cite{boroujeni2017flexible}. None of the above mentioned A* based motion planners does not consider energy consumption. 

An application of A*-search for finding energy optimal velocity trajectory for an electric bicycle was presented in \cite{grossoleil2012practical}. Authors used kinetic and potential energy as well as a rolling resistance, but for estimation of air drag resistance, the authors used an upper instead of a lower bound. The results were suboptimal with a difference of around 1.2\% from another optimal control strategy. In \cite{chevrant2014search} authors introduced air drag and a time proportional cost without interdependence. For air drag, they used a tunable minimum velocity, and for the time-proportional cost, they computed a minimum time based on the maximum velocity. This approach leads to a not so precise estimation and a loss of admissibility if velocities are lower than tuned value. Another use of A*-search was presented in \cite{flehmig2015energy} but authors did not reveal the computation of their heuristics. From the results, it is clear that some of the heuristics are giving a suboptimal solution. 

In this work, a novel heuristic function for A*-based vehicle velocity trajectory planning is introduced. The optimal cruising velocity, which minimizes air drag and time-proportional cost influence, is defined. The resulting, optimal cruising velocity is then used as a basis for the introduced heuristic function.


\section{Problem definition}

The problem of finding an energy efficient vehicle velocity trajectory can be considered as an optimal motion planning problem. This problem has differential constraints defined by a vehicle dynamic model, state-dependent constraints which have to be satisfied and a cost function which has to be minimized.

\subsection{System model}

Since low model complexity is crucial for efficient optimization, the vehicle is modeled as a lumped mass. The vehicle is represented by two states: $s$ represents traveled longitudinal distance and $v$ represents longitudinal velocity of the vehicle in the lane and vehicle dynamics as in (\ref{eq:veh1})-(\ref{eq:veh3}). 
\begin{align}\label{eq:veh1}
 \dot{s}(t) =& v(t),\\
 \dot{v}(t) =& \frac{F_m(t) - F_r(t)}{m},\\ \label{eq:veh3}
 F_r(t) =& \frac {1}{2}\rho _a c_d A_f v(t)^2 + c_r m g \cos \left( \alpha \left(s(t) \right) \right)\hdots\\\notag &\quad \hdots+ m g \sin \left( \alpha \left(s(t) \right) \right).
\end{align}
The propulsion force of the electric motor is denoted by $F_m$ and resistance forces by $F_r$. As it is shown in (\ref{eq:veh3}), resistance forces considered include air drag resistance, rolling resistance and road slope ($\alpha$) related gravity force. The propulsion element, an electric motor, with inner torque $T_m$ is modeled statically by:
\begin{equation} \label{eq:torque}
F_m(t)=\frac{kT_m(t)\eta^{\mathit{sign}(T_m(t))}}{r_W},\ \text{and} \ k = \frac{2 r_W \pi \omega_m(t)}{v(t)},
\end{equation}
where $\eta$ is an efficiency coefficient, $k$ is a combined transmission ratio of the powertrain, $r_W$ is the radius of the wheels and $\omega_m$ is the rotational speed of the motor.

\subsection{Cost function}

The cost function has to reflect the initial requirement of minimal energy consumption. When considering only propulsion power in the cost function, energy-efficient behavior results in a zero velocity trajectory. To avoid this, some authors introduced a term to the cost function to weight the traveling time (\cite{sciarretta2015optimal}). The weighting coefficient is then tuned such that the travel time is comparable to times achieved by human drivers. However, energy consumption includes energy used for auxiliaries $P_{aux}$ (e.g. infotainment, air conditioning, etc.) which can be approximated with a constant load. This load corresponds to a weighting coefficient presented in other works, but it is determined experimentally and not tuned.
Consequently, the energy used is the integral of the power over time of the trip represented by:
\begin{equation} \label{eq:cost_time}
E_\mathrm{min} = \int_{t_i}^{t_f} \left( \omega(t)T_m(t)+P_{aux} \right)dt.
\end{equation}

Instead of using time as a variable for integration, the distance can be used too (\cite{saerens2012optimal}). This offers some advantages for solving as final time is not known, and final distance is, and road slope appears as a function of distance.
\begin{equation} \label{eq:cost_dist}
E_\mathrm{min} = \int_{s_i}^{s_f} \left( \frac{k T_m(s)}{2 r_W \pi}+\frac{P_{aux}}{v(s)} \right)ds.
\end{equation}

\subsection{Complexity}
To solve the optimization task numerically, using graph searching methods, state discretization is necessary. By increasing the number of system states considered in the optimization problem, complexity increases rapidly. Additionally, and even more problematic, the number of possible transitions that need to be evaluated in each step increases significantly.  Bellman called this problem the "curse of dimensionality" (\cite{bellman1954theory}). 

\section{OPTIMAL MOTION PLANNER}

In this work two different motion planners are used, one based on DP as discussed in  \cite{myBookCh2017} and one using A* which is introduced in this paper.

\subsection{A* search}

A* is one of the earliest yet one of the most used methods for path planning (\cite{hart1968formal}). It is based on the well-known, Dijkstra's algorithm but uses heuristics to guide exploration to the nodes which (according to heuristic) lead to the solution quicker and therefore has a better performance compared do Dijkstra's algorithm. A* is an optimal method for finding the optimal path for a certain admissible heuristic function.

Starting from the \textit{initial node}, which is chosen as the initial \textit{current node}, all neighbors are determined and added to the $\OPEN$ list. From the $\OPEN$ list, the node with the lowest cost is chosen to be the next \textit{current node}. This is repeated until the \textit{goal node} is reached or the whole graph is explored. The cost $f(n)$ is computed using (\ref{eq:f}) where $g(n)$ represents the cost to travel from \textit{initial node} to \textit{current node} (\textit{cost-to-come}) and $h(n)$ the cost to travel from \textit{current node} to \textit{goal node} (\textit{cost-to-go}), estimated using some heuristic function.
\begin{equation} \label{eq:f}
f(n)=g(n)+h(n).
\end{equation}
There can be several issues when implementing A*. To accelerate accessing the minimum of all open nodes, the $\OPEN$ list is usually organized as some kind of a priority queue. In this work, the $\OPEN$ list is implemented as a binary heap with a hash table. A binary heap is a structure which keeps the minimum element at the top and therefore it is easy to access when selecting the next \textit{current node}. The hash table keeps a position of every node on the heap for an easier accessing when comparing with newly expanded nodes. Besides, the hash table keeps information if a node is already in the $\CLOSED$ list.

\subsection{Cost-to-come} 
Cost-to-come $g(n)$ is the exact cost necessary to come from the \textit{initial node} to the \textit{current node} $n$. It is usually computed cumulatively. One step transition cost is added to the \textit{parent node}'s \textit{cost-to-come}. To compute transition costs the vehicle model (\ref{eq:veh1})-(\ref{eq:torque}) is used along with the cost function (\ref{eq:cost_dist}).

\subsection{Search space}
As the planning task is solved by A*, graph searching method, state discretization is necessary and an appropriate graph has to be constructed. The graph is constructed by discretizing distance ($s$) and velocity ($v$) states. Instead of distance, time can be considered also, but choosing distance is useful when the final time is not fixed. This choice brings several disadvantages e.g. if a trajectory passes zero velocity, it is impossible to compute time spent in that state. By assuming positive velocities only, the graph is directed. To consider dynamic obstacles, an additional dimension for time $t$ should be used (\cite{fraichard1993dynamic}), but this is not in the focus of this work. 
\begin{figure}
  \includegraphics[angle =-90, trim=97 130 230 105,clip, width=\linewidth ]{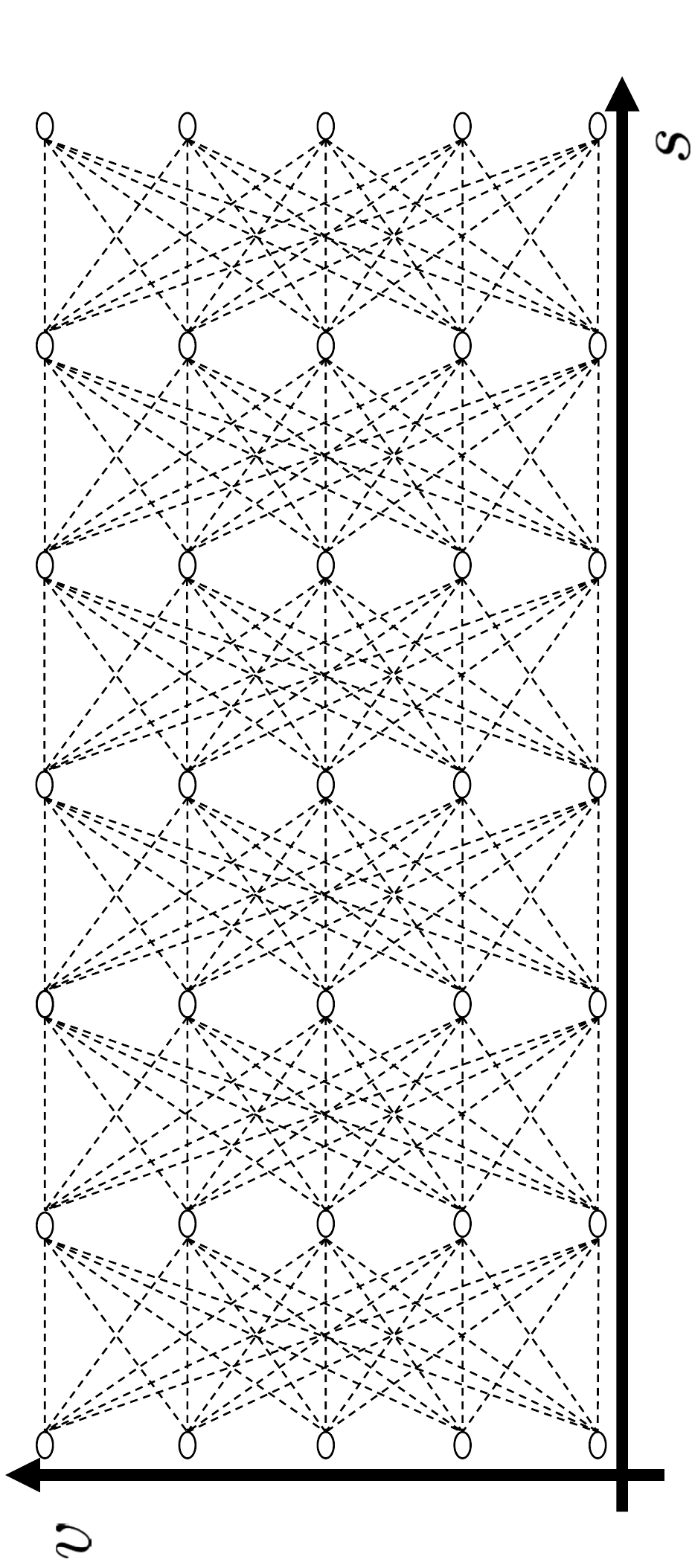}
  \caption{Search space constructed as a graph.}
  \label{fig:graph}
\end{figure}

\subsection{Heuristic function for cost-to-go estimation}

The heuristic function is used to estimate the cost needed to travel from any node (point in state space), defined by initial velocity and position ($v_i$ and $s_i$) to the goal node. As it is shown in \cite{hart1968formal}, if the heuristic function is admissible (underestimates the \textit{cost-to-go}), the result of A* search is the optimal trajectory. For the shortest path search, the usual heuristic function is the Euclidian distance. To find the energy optimal velocity trajectory, the heuristic function must always underestimate the energy needed to drive from any node to the goal node. The precision of the heuristic influences the size of the explored space and therefore the efficiency of the search. In general, as a heuristic function tends to the exact \textit{cost-to-go}, the explored space shrinks and the search time shortens. In the ideal case, if the heuristic function gives the exact cost-to-go, only nodes which belong to the optimal trajectory would be explored.

\section{HEURISTIC FUNCTIONS}
The energy needed to drive from the initial state ($v_i$ and $s_i$) to the final state ($v_f$ and $s_f$) based on model (\ref{eq:veh1})-(\ref{eq:veh3}) and the cost function (\ref{eq:cost_dist}), can be presented as the sum of the work needed for acceleration, overcoming air drag, rolling resistance,  road slope and the energy used by the auxiliaries (or operational costs). 
\begin{equation} \label{eq:cost_time}
W = W_{acc}+W_\alpha+W_{roll}+W_{air}+E_{aux}.
\end{equation}
Assuming no energy is needed to drive to the final state would provide admissible heuristic if recuperation would not be allowed, as the work is always positive. Since recuperation is allowed, under certain conditions, the work can be negative too. So, assuming a value of zero could lead to a suboptimal solution. Air drag, rolling resistance, and auxiliary-related works are always positive and therefore can be underestimated by zero. Nevertheless, this is an imprecise estimation, far from the real value, leading to the loss of A* efficiency as informed search.

\subsection{Heuristic function for acceleration related work}
If no efficiency losses are considered, the work needed to accelerate the vehicle moving in the horizontal plane results in a change of kinetic energy. So the work needed for acceleration can be described only by the initial and the final velocity states, $v_i$ and $v_f$ as in (\ref{eq:heur_acc}). If losses are present, the needed energy will only be higher.
\begin{equation} \label{eq:heur_acc}
W_{acc} \geq \Delta E_k = m \frac{v_f^2-v_i ^2}{2}.
\end{equation}
If the vehicle is driving on the horizontal plane, without slopes, this work alone expresses admissible heuristic. If this is not the case, the influence of road slope has to be included additionally.

\subsection{Heuristic function for slope related work}
Overcoming road slope results in climbing up or down a hill and therefore increases or decreases the potential energy of the vehicle. In this case, the work can be described as:
\begin{equation} \label{eq:cost_time}
W_\alpha+W_{acc} \geq \Delta E_k + \Delta E_p= \Delta E_k +mg(h(s_f)-h(s_i))
\end{equation}
Where $h(s)$ is the road elevation profile dependent on the distance traveled ($s$).

\subsection{Heuristic function for rolling resistance related work}
As the rolling resistance force is modeled as a constant, it does not depend on vehicle driving trajectory. It can be exactly computed for the specific road segment and included in the heuristics. Assuming the rolling resistance work as a constant provides a more precise, and still admissible heuristic, but cannot be used alone, apart from $\Delta E_k + \Delta E_p$.
\begin{equation}
W_{roll} = c_r m g \int_{s_i}^{s_f} \cos \left( \alpha \left(s \right) \right)ds
\end{equation}

\subsection{Heuristic function for air-drag and auxiliary power}\label{sec:wa}
\begin{figure}[b]
  \includegraphics[trim=50 6 40 15,clip, width=\linewidth]{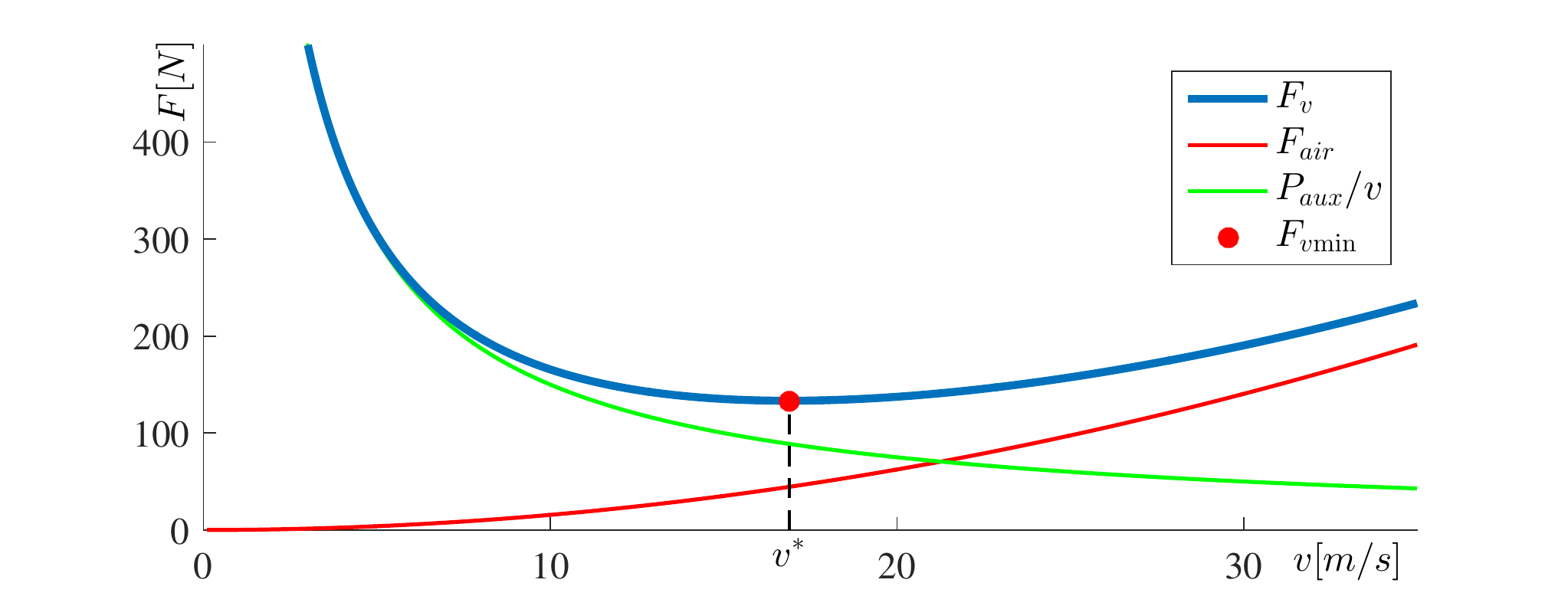}
  \caption{Velocity-dependent virtual force $F_v$, generated by combined air drag resistance and auxiliary power.}
  \label{fig:conv}
\end{figure}
Above mentioned heuristic functions are used for constructing the state-of-the-art admissible heuristic function. So far, based on the authors best knowledge, none of the heuristics successfully incorporated air-drag resistance and auxiliary power to provide admissible heuristic.

As opposed to $W_{acc}$, $W_\alpha$ and $W_{roll}$, which can be expressed by the initial and the final states only, $W_{air}$ and $E_{aux}$ depend on velocity trajectory between the initial and the final state.
\begin{align} \label{eq:cost_time}
W_{air} + E_{aux} =& \int_{s_i}^{s_f} \left(\frac {1}{2}\rho _a c_d A_f v(s)^2 + \frac{P_{aux}}{v(s)}\right)ds,\\
W_{air} + E_{aux} =& \int_{s_i}^{s_f} F_v \left(v(s)\right)ds,\\
F_v =& \frac {1}{2}\rho _a c_d A_f v(s)^2 + \frac{P_{aux}}{v(s)}.
\end{align}
$F_v$ can be considered as a virtual force consisting of air drag force and virtual force by $P_{aux}$. $P_{aux}$ can include also some economic costs such as hourly rate ($hr$) of the operator, vehicle renting, etc. These are then divided by the current electricity price ($ep$) to get the power equivalent:
\begin{equation} \label{eq:cost_time}
P_{tot} = P_{aux}+\frac{hr}{ep}.
\end{equation}
As it can be seen in Fig. \ref{fig:conv} based on the convex shape of $F_v$, we can conclude that there exists one optimal velocity for which $F_v$ is minimized, noted by $v^*$. 

\subsubsection{Optimal cruising velocity.} This optimal cruising velocity $v^*$ minimizes the virtual force $F_v$: 
\begin{equation}
F_{v\mathrm{min}} = F_v(v^*).
\end{equation}
To find the $v^*$, the derivative of $F_v$ at $v^*$ should be zero:
\begin{equation}
\frac{\partial F_v(v^*)}{\partial v} = 0 \quad \Rightarrow \quad v^* = \sqrt[3]{\frac{P_{aux}}{\rho _a c_d A_f }}.
\end{equation}
The optimal cruising velocity is unique for the vehicle (and it's hourly rate $hr$). It depends on vehicle aerodynamic shape, air density and $P_{aux}$ (eventually, motor efficiency can be included). 
The minimum virtual force $F_{v\mathrm{min}}$, can be used to estimate the lower bound of $W_{air} + E_{aux}$, defined only with the initial and the final states. 
\begin{equation}\label{eq:conv}
W_{air} + E_{aux} \geq W_A = \int_{s_i}^{s_f} F_{v\mathrm{min}}ds = F_{v\mathrm{min}}(s_f-s_i).
\end{equation}
As it was shown in (\ref{eq:conv}), air drag and auxiliary-related work, $W_{air} + E_{aux}$ for any trajectory, is always greater or equal to the work when driving with constant velocity $v^*$ (Fig. \ref{fig:conv}). This work is computed as a product of the minimum virtual force $F_{v\mathrm{min}}$ and the distance traveled.
\begin{figure}
  \includegraphics[angle =-90, trim=95 110 223 100,clip, width=\linewidth ]{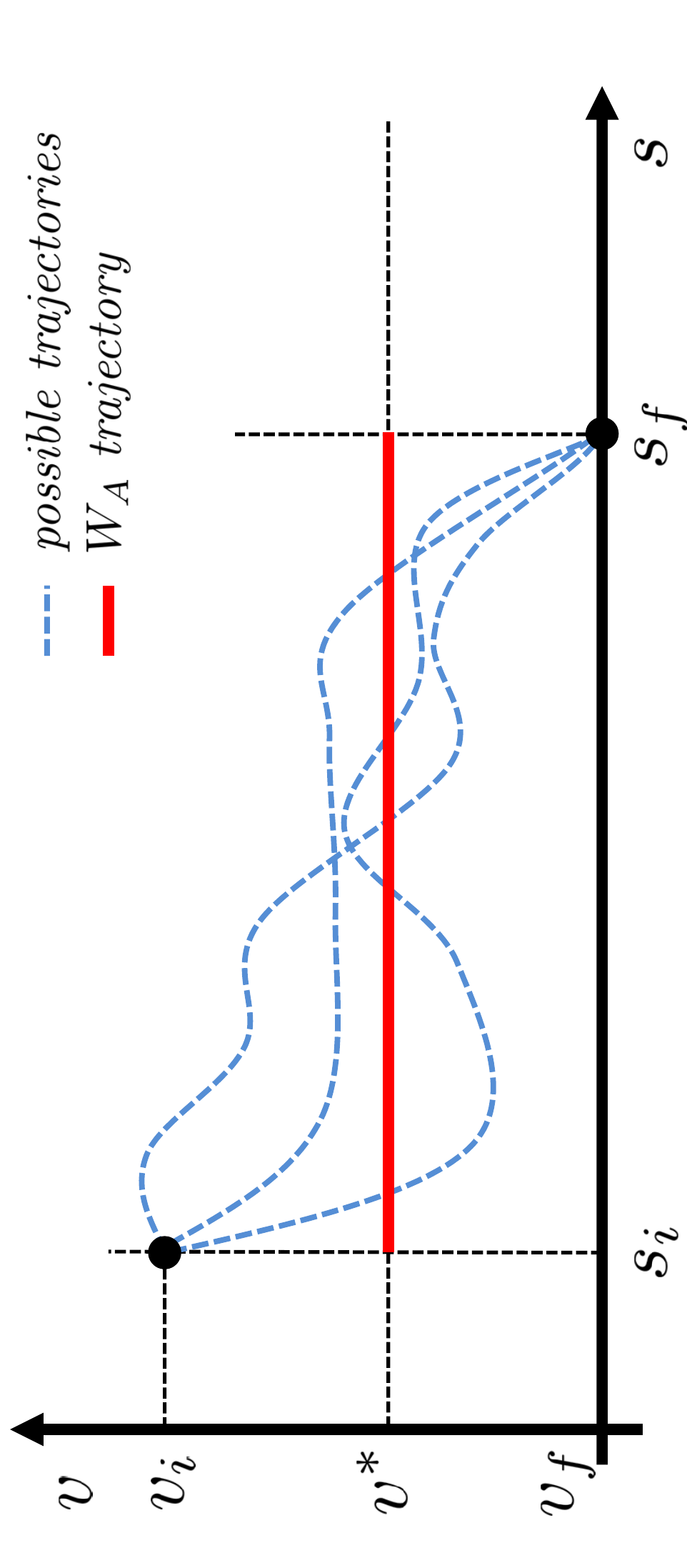}
  \caption{Some of possible trajectories for driving from initial to final state, and lower bound trajectory for $W_A$ computation.}
  \label{fig:fig3}
\end{figure}

A tendency to drive with optimal cruising velocity can be noticed in \cite{xu2017instantaneous} where authors analyzed a conventional vehicle using approximated engine fuel injection rate map. As in conventional vehicles auxiliary power comes from the alternator, which represents an additional load on the engine, the results are similar. In our work, a theoretical explanation of this phenomenon (optimal cruising velocity), reasoning and ifluencing factors and relations are provided.
\subsection{Improved heuristic function for air-drag and auxiliary power} \label{sec:wai}
\begin{figure}[t]
  \includegraphics[angle =-90, trim=95 110 215 100,clip, width=\linewidth ]{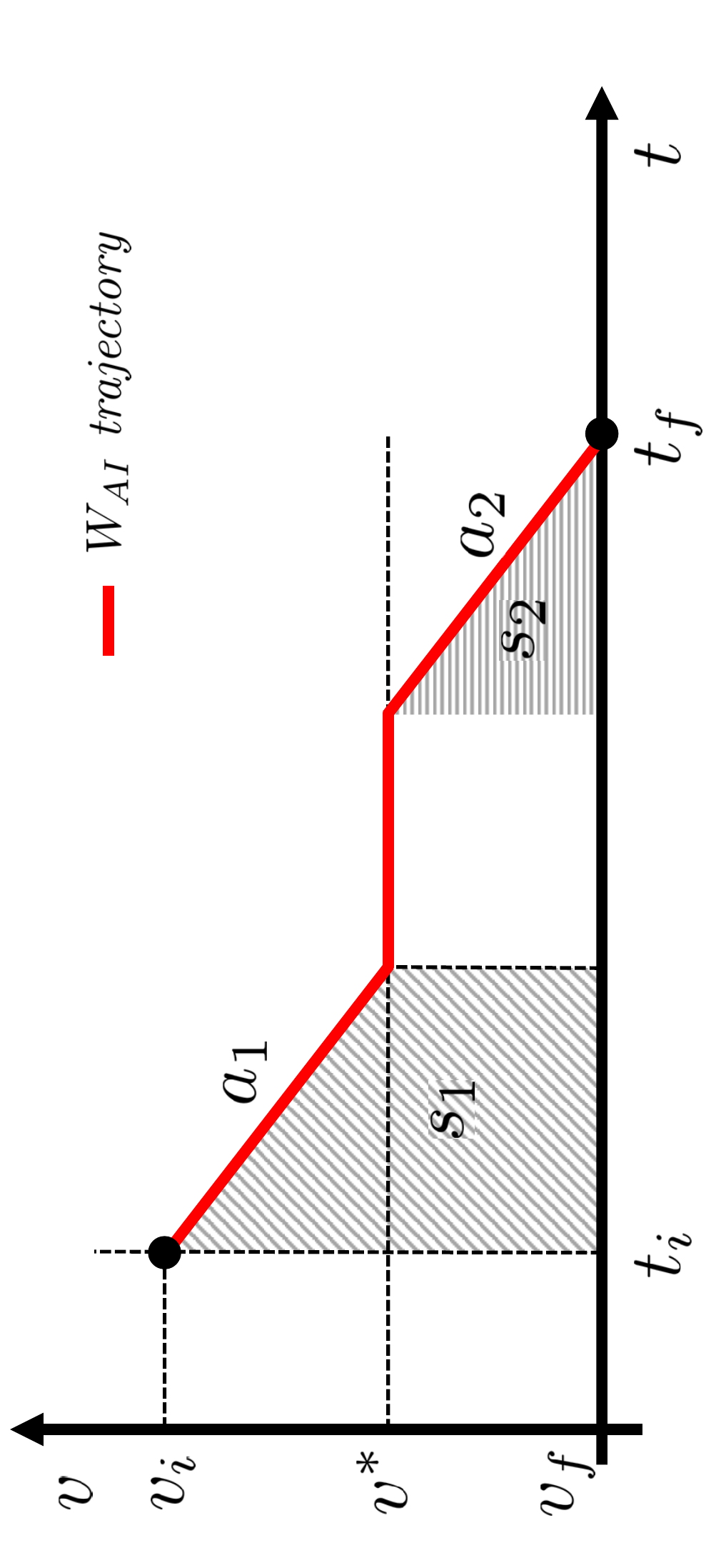}
  \caption{Velocity trajectory for $W_{AI}$ computation.}
  \label{fig:fig4}
\end{figure}

For situations where the initial or the final velocity is not equal to the optimal cruising velocity $v^*$, the air drag influence can be more precisely estimated by including acceleration periods to reach $v^*$. As acceleration is limited, this transition is not instantaneous. It is important to note that work needed for acceleration is not computed in this trajectory, as it is included in kinetic and potential energy. Only energy $W_{air} + E_{aux}$, as if the vehicle would move on this trajectory is considered. A general trajectory with acceleration and deceleration phases is depicted on picture Fig. \ref{fig:fig4}. 
For simpler visualization, $t$ is used as the axis instead of $s$.

While driving with uniform acceleration, from some velocity $v_1$ to some velocity $v_2$, the distance traveled can be computed as:
\begin{equation}\label{eq:dist}
s=\frac{v_2 ^2 - v_1 ^2}{2 \cdot a}.
\end{equation}
For computing $s_1$ form Fig. \ref{fig:fig4}. general formula (\ref{eq:dist}) is used, with $v_1 = v_i$, $v_2 = v^*$ and $a = a_1$. For computing $s_2$ it is with $v_1 = v^*$, $v_2 = v_f$, and $a = a_2$. Depending on the initial and the final velocity accelerations $a_1$ and $a_2$ are determined as: 
\begin{equation}
\resizebox{.87\hsize}{!}{
$a_1 = \begin{cases}
			a_{min} , v_i > v^*,\\ 
			a_{max} , v_i < v^*,
\end{cases}\ \text{and} \quad
a_2 = \begin{cases}
			a_{min} , v_f < v^*,\\ 
			a_{max} , v_f > v^*,
\end{cases}$}
\end{equation}
The work $W_{AT} = W_{air} + E_{aux}$, when driving with uniform acceleration, for time $T$, from velocity $v_1$ to $v_2$:
\begin{equation} \label{eq:wat1}
W_{AT} = \int_0^{T} \left(\frac {1}{2}\rho _a c_d A_f v(t)^3 + P_{aux}\right)dt,\\
\end{equation}
\begin{equation}\label{eq:wat2}
 \text{with} \ v(t) = v_1+at, \ \text{and} \ T = \frac{v_2-v_1}{a}, 
\end{equation}
\begin{equation}\label{eq:wat3}
\resizebox{.87\hsize}{!}{
$W_{AT} = \frac {1}{2}\rho _a c_d A_f \left(v_1^3T + \frac{3v_1^2aT^2}{2}+\frac{a^3T^4}{4} \right) + P_{aux}T,$}
\end{equation}

Using relations (\ref{eq:wat2}) and (\ref{eq:wat3}), the cost of driving with uniform acceleration between two velocities can be computed. The total work $W_{AI}$, in this case, can be computed as the sum of two works (acceleration and deceleration) and the driving with constant velocity $v^*$: 
\begin{equation}\label{eq:WAI}
W_{AI} = W_{AT1} + F_{v\mathrm{min}} \cdot (s_f - s_i - s_1 - s_2)+W_{AT2}.
\end{equation}

\begin{figure}[t]
 \includegraphics[angle =-90, trim=95 110 215 100,clip, width=\linewidth ]{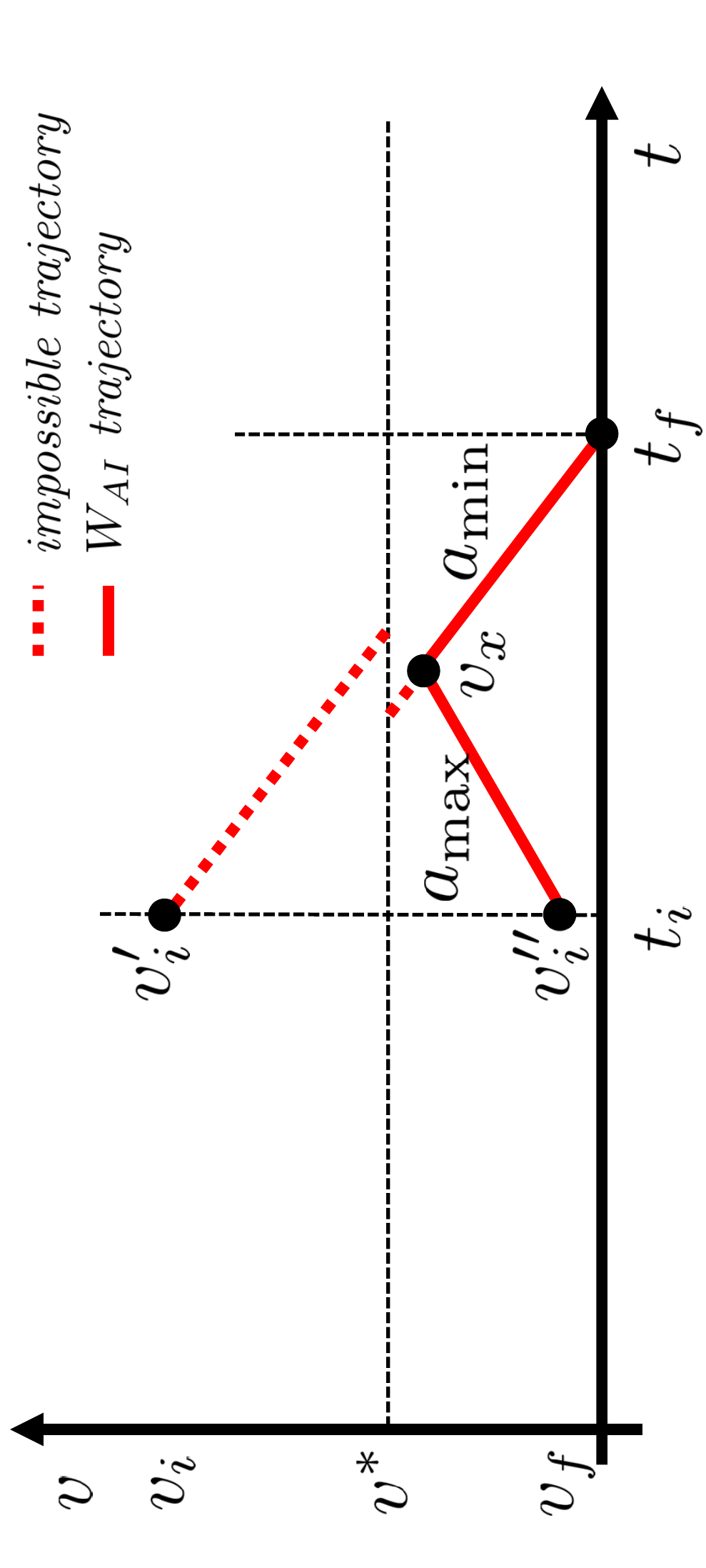}
  \caption{Special cases of velocity trajectories for $W_{AI}$ computation.}
  \label{fig:fig5}
\end{figure}

If $s_1+s_2 > s_f - s_i$, the optimal velocity cannot be reached as shown in Fig. \ref{fig:fig5}. There are two subcases in this case. 1) When $v_i$ and $v_f$ are on opposite sides of $v^*$. This is a trivial case; as the final state is not reachable for limited acceleration and therefore the cost is infinite. This case is marked with dashed red line on Fig. \ref{fig:fig5}. 2) When both $v_i$ and $v_f$ are either smaller or greater than $v^*$. There will be no constant velocity driving part as $v^*$ will not be reached, it will accelerate until it reaches $v_x$ and then decelerate, or the opposite. This is depicted with a solid red line on Fig. \ref{fig:fig5}.
By using the $s_1+s_2= s_f - s_i$, the velocity at which acceleration changes sign can be determined as:
\begin{equation}\label{eq:vx}
v_x = \sqrt{\frac{2a_1a_2s+a_2v_i^2-a_1v_f^2}{a_2-a_1}}.
\end{equation}
Velocity $v_x$ is then used instead of $v ^*$ for computing $W_{AT1}$ and $W_{AT2}$. The work $W_{AI}$, in this case, can be computed as the sum of two works (acceleration and deceleration):
\begin{equation}\label{eq:WAIx}
W_{AI} = W_{AT1} + W_{AT2}.
\end{equation}

\subsection{Motor efficiency}
The motor efficiency certainly influences the energy used. It does not have to be considered to have admissible heuristics, but it improves precision and therefore search efficiency. Applying the maximum motor efficiency $\eta_\mathrm{max}$ as a constant efficiency for any operating point will result in underestimating heuristic which is more precise than when omitting efficiency. The motor efficiency can be applied only to power which flows through the motor. Therefore, it should not be applied to $P_{aux}$ as this energy does not flow through the motor.
\subsection{Conclusion on heuristics}
Combining all relevant work elaborated in the introduction the best state-of-the-art, admissible heuristic, can be constructed as a sum of kinetic and potential energy, rolling resistance loss with applied motor efficiency as in (\ref{eq:H_SoA}) - (\ref{eq:H_SoA2}). 
\begin{align}\label{eq:H_SoA}
W_{tot} =& \Delta E_k + \Delta E_p + W_{roll},\\ \label{eq:H_SoA2}
h_{SoA} =& W_{tot} \cdot \eta_{max}^{-\mathit{sign}(W_{tot})}.
\end{align}
The proposed heuristic, improves this state-of-the-art heuristic by additional work for air drag resistance and auxiliary energy (and operation costs) as shown in the sections \ref{sec:wa} and \ref{sec:wai}. The efficiency cannot be applied to the proposed heuristic as auxiliary power does not flow trough motor. The proposed heuristic is computed as: 
\begin{equation}\label{eq:H_pro}
h_{pro} = h_{SoA} +W_{AI}.
\end{equation}

\section{SIMULATION RESULTS}

\begin{figure}
  \includegraphics[angle =0, trim=35 2 50 20,clip, width=\linewidth ]{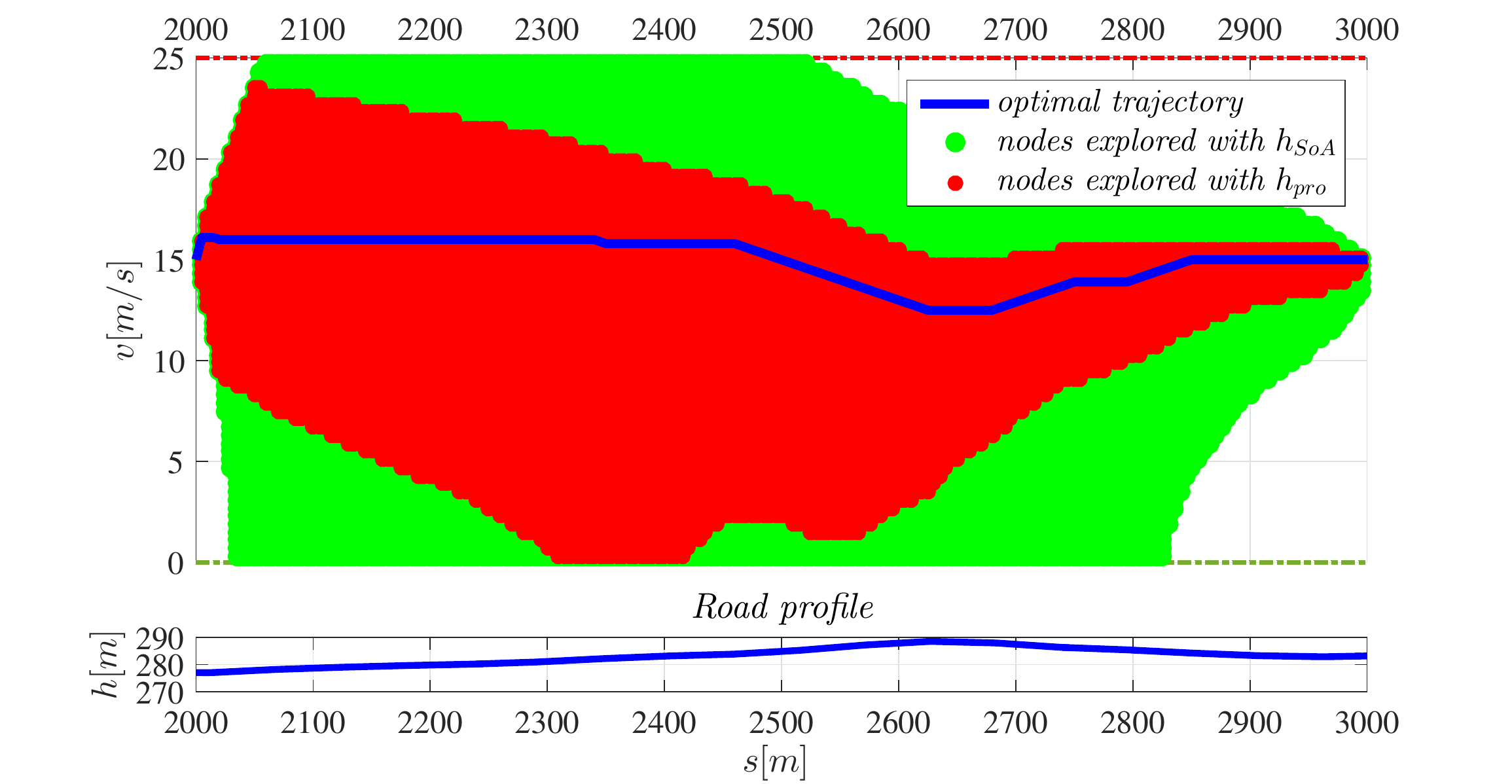}
  \caption{Space searched by using different heuristics.}
  \label{fig:sim}
\end{figure}

For a simulation validation,three methods are compared: DP, A*-search with state-of-the-art heuristic and A*-search with the proposed heuristic. Methods are used to plan the optimal velocity trajectory for driving on a $1 km$ long segment of the A9 highway, in the vicinity of Graz, Austria. All methods result in the exact same optimal velocity trajectory (Fig. \ref{fig:sim}), which verifies implementations of both  algorithms (DP and A*), and admissibility of both heuristics. Total energy needed for this trip was computed as \textit{422,8 kJ}.
A comparison between state-of-the-art heuristic and proposed heuristic, is made based on the error of the estimation and search efficiency represented by the number of explored nodes. DP provides the exact \textit{cost-to-go} computation which is used as a reference for computation of the estimation error of the heuristics. 
\begin{table}[h]
\begin{center}
\caption{Comparison of the simulation results.}\label{tb:table}
\begin{tabular}{cccc}
  & DP & A* with $h_{SoA}$ & A* with $h_{pro}$ \\\hline
Nodes explored & 50200 & 41125 & 25052 \\
Average error [kJ] & 0 & -84,2 & -15,2 \\ 
Minimum error [kJ] & 0 & -173,5 & -37,6 \\ 
Maximum error [kJ] & 0 & 0 & 0 \\\hline
\end{tabular}
\end{center}
\end{table}
The results are shown in Table \ref{tb:table} along with the results from DP. The first column contains results of using DP. As DP provides the exact \textit{cost-to-go} errors are zero. The second and the third column contain results of A*-search with state-of-the-art heuristics (\ref{eq:H_SoA2}) and the proposed heuristics (\ref{eq:H_pro}) respectively.  As it is shown in Table \ref{tb:table}., the average error improved almost \textbf{5 times} and number of explored nodes almost \textbf{2 times} when the proposed heuristic ($h_{pro}$) is used, compared to the state-of-the-art heuristic ($h_{SoA}$). As it can be seen, the error is always negative, which means that the heuristics always underestimates the \textit{cost-to-go}. A*, in general, examines a smaller number of nodes than DP, and it further decreases as the precision of the heuristics increases. Figure \ref{fig:sim} presents the optimal trajectory and space searched by state-of-the-art and proposed heuristics from the same simulation. It is important to note that, the number of explored nodes for DP depends on the maximum velocity used; and for A* it does not.

\section{Conclusion}
The proposed heuristic showed to effectively incorporate air drag and auxiliary power (and operational costs) for energy(economic)-optimal motion planning. It was shown that proposed heuristics significantly improves the precision of the estimation, which decreases the number of examined nodes compared to the state-of-the-art heuristic. Moreover, the proposed heuristic explores significantly less trajectories compared to dynamic programming approach.
Beside motion planning, model-based heuristics could be used for the estimation of the range of electric vehicles. As the proposed heuristic always underestimates the cost-to-go, it could be also used to eliminate overestimating of the vehicle remaining range. Additionally, the optimal cruising velocity, could be used to advise drivers to drive more efficiently.

\begin{ack}
\begin{footnotesize}
The project leading to this study has received funding from the European Union's Horizon 2020 research and innovation programme under the Marie Skłodowska-Curie grant agreement No 675999, ITEAM project.\par
VIRTUAL VEHICLE Research Center is funded within the COMET – Competence Centers for Excellent Technologies – programme by the Austrian Federal Ministry for Transport, Innovation and Technology (BMVIT), the Federal Ministry of Science, Research and Economy (BMWFW), the Austrian Research Promotion Agency (FFG), the province of Styria and the Styrian Business Promotion Agency (SFG). The COMET programme is administrated by FFG.\par
\end{footnotesize}
\end{ack}

\bibliography{ifacconf}             
                                                  
\end{document}

%% file: macros/symbols.tex
\newcommand{\OPEN}{\textsc{Open}}
\newcommand{\CLOSED}{\textsc{Closed}}	